%
%

%
%


\magnification 1200
\input amstex
\documentstyle{amsppt}
\NoBlackBoxes
\NoRunningHeads

\vsize = 9.4 truein

\define\vre{\varepsilon}

\define\hs{homogeneous space}
\define\df{\overset\text{def}\to=}
\define\un#1#2{\underset\text{#1}\to#2}
\define\br{\Bbb R}
\define\bn{\Bbb N}
\define\bz{\Bbb Z}
\define\bq{\Bbb Q}
\define\bc{\Bbb C}
\define\bk{\Bbb K}

\define\da{Diophantine approximation}
\define\de{Diophantine exponent}

\define\ve{\bold e}
\define\vx{\bold x}
\define\vy{\bold y}
\define\vz{\bold z}
\define\vu{\bold u}
\define\vv{\bold v}

\define\vs{\bold s}
\define\vw{\bold w}

\define\vq{\bold q}

\define\vc{\bold c}
\define\vf{\bold f}
\define\vg{\bold g}
\define\vh{\bold h}
\define\vt{\bold t}

\define\spr{Sprind\v zuk}

\define\nz{\smallsetminus \{0\}}

\define\cag{$(C,\alpha)$-good}

\define\kgt{Khintchine-Groshev Theorem}

\define\GL{\operatorname{GL}}
\define\SL{\operatorname{SL}}
\define\im{\operatorname{Im}}

\topmatter
\title 
Baker-Sprind\v zuk conjectures for complex analytic manifolds 
\endtitle  

\author { Dmitry Kleinbock} \\ 
  { \rm 
   Brandeis University} 
\endauthor

    \address{ Dmitry Kleinbock,  Department of
Mathematics, Brandeis University, Waltham, MA 02454-9110}
  \endaddress

\email kleinboc\@brandeis.edu \endemail

  \thanks  
Supported in part by NSF
Grant DMS-0072565.
\endthanks

\abstract 
We show a large 
class of analytic submanifolds of $\bc^n$ to be strongly extremal. This generalizes
V.~\spr's solution of the complex case of Mahler's Problem, and settles
complex analogues of conjectures made in the 1970s by Baker and \spr.
The  proof is based on 
a variation of quantitative nondivergence estimates for
quasi-polynomial
flows on the space
of lattices.

 \endabstract


\date October 16, 2002
\enddate

\endtopmatter
\document

\heading{1. Introduction}
\endheading 

The circle of problems that the present paper belongs to dates back to the 1930s,
namely, to K.~Mahler's work on a classification of transcendental real and complex
numbers. For a polynomial $P(x) = a_0 + a_1x + \cdot + a_nx^n\in\bz[x]$, let us denote
by $h_P$ the {\sl height\/} of $P$, that is, $h_P \df\max_{i =
0,\dots,n}|a_i|$.  It can be easily shown using Dirichlet's Principle that for any
$z\in \bc$ and any $n\in\bn$ there exists a positive constant
$c(n,z)$ such that
$$
|P(z)| < c(n,z)h_P^{-v} \quad\text{for infinitely many
}P\in\bz[\cdot]\quad\text{with }\deg P\le n\,,\tag 1.1
$$
where $v = \frac{n-1}2$, and one can take $v = n$ if $z\in \br$. Mahler's
Conjecture \cite{M}, proved in 1964 by V.~\spr\ \cite{S1, S2, S3}, states that for almost
every
$z\in\bc$ (resp.~$z\in\br$), the values $v = \frac{n-1}2$
(resp.~$v = n$) in (1.1) cannot be increased. Loosely put, this result show
that almost all real/complex numbers are `as far from being algebraic as they could
possibly be'.

Let us restate the aforementioned results by defining the {\sl \de\/} $\omega(\vz)$ of
$\vz\in\bc^n$ by
$$
\omega(\vz) \df \sup\left\{v > 0\left|\aligned
 |\vz\cdot\vq + p|   \le \|\vq\|^{-v}&\text{ \ for}\\\text{infinitely many }\vq\in
\bz^n&,\  p\in\bz\endaligned\right.\right\}\,,\tag 1.2
$$
where  $\|\vq\|$ stands for
$\max_{i}|q_i|$. Dirichlet's Principle  and the Borel-Cantelli Lemma imply that
$\omega(\vz)$ is:

\roster
\item"$\bullet$" not less than $\frac{n-1}2$ for all $\vz\in\bc^n$, and equal to
$\frac{n-1}2$ for almost all $\vz\in\bc^n$; 
\item"$\bullet$" not less than $n$ for all $\vz\in\br^n$, and equal to
$n$ for almost all $\vz\in\br^n$.
\endroster

Let us use the notation $\bk$ for $\br$ or $\bc$. Following a terminology introduced by
\spr, say that a map $\vf$ from an open subset $U$ of $\bk^d$ to $\bk^n$ is {\sl
extremal\/}  if the \de\ of $\vf(\vx)$ is for almost every $\vx\in U$ equal to that of a
generic point of $\bk^n$ (that is, to $n$ if $\bk = \br$ or to $\frac{n-1}2$ if $\bk =
\bc$), and that a smooth submanifold of $\bk^n$ is {\sl extremal\/} if so are all its
parametrizing maps. \spr's result therefore states that 
the  curve  
 $$
\Cal M = \{(z,z^2,\dots,z^n)\mid z\in \bk\}\subset \bk^n\tag 1.3
$$ 
is extremal. Thus a natural generalization of Mahler's Problem is to look for
general conditions sufficient for the extremality of a manifold/map. 

\smallskip

Another extension arises if one replaces  $\|\vq\|$ in (1.2) by the `geometric
mean' of the components of $\vq$. More precisely, define
the {\sl multiplicative \de\/} $\omega^\times(\vz)$ of
$\vz\in\bc^n$ by
$$
\omega^\times(\vz) \df \sup\left\{v > 0\left|\aligned
 |\vz\cdot\vq + p|   \le \big(\Pi_+(\vq)\big)^{-v/n}\text{ \ for}\\ \text{infinitely
many }\vq\in
\bz^n,\  p\in\bz\endaligned\right.\right\}\,,
$$
where 
$$
\Pi_{\sssize +}(\vq) 
= \prod_{i = 1}^n |q_i|_{\sssize +}\quad\text{and}\quad|q|_{\sssize +}=\max(|q|,
1)\,.
$$
Since $
\Pi_{\sssize +}(\vq)$ is not greater than $\|\vq\|^n$ for any nonzero $\vq$, one 
has $
\omega^\times(\vz) \ge 
\omega(\vz)$ for any $\vz\in\bc^n$. However, one can show by a Borel-Cantelli-type
argument that the two exponents agree for almost all $\vz$; more precisely, that $
\omega^\times(\vz)$ is equal to $\frac{n-1}2$ for almost all $\vz\in\bc^n$ and 
to
$n$ for almost all $\vz\in\br^n$.

The multiplicative analogue of the notion of extremality is usually referred to as
strong extremality. Let us say that a map $\vf:U\to 
\bk^n$ is {\sl strongly
extremal\/}  if  $
\omega^\times\big(\vf(\vx)\big)$ is for almost every $\vx\in U$ equal to $\frac{n-1}2$
or $n$ for $\bk = \bc$ of $\br$ respectively, with a similar definition for strongly
extremal 
submanifolds of $\bk^n$.

\smallskip

Identifying extremal and strongly
extremal manifolds has been one of the central issues of 
metric \da\ for the
last 40 years. However most of the activity revolved around the case $\bk = \br$.
In his 1980 survey of the field \cite{S4}, Sprind\v zuk conjectured that a
real analytic submanifold $\Cal M$ of $\br^n$ is strongly extremal whenever it  is not
contained in any proper affine subspace of $\br^n$; or, equivalently, that a real
analytic map
$\vf = (f_1,\dots,f_n)$ is strongly extremal if 
$$
1,f_1,\dots,f_n\text{ are linearly independent
over }\br\,.\tag 1.4
$$ 
The latter condition, loosely put, says that $\Cal M$ `remembers' 
the dimension of the space it is imbedded into, and the conjecture asserted that $\Cal M$
must also `remember' the \de\ of a generic point of $\br^n$. The special case of $\Cal M$
of the form (1.3) and $\bk = \br$, that is, the multiplicative analogue of Mahler's
original problem, was conjectured earlier by A.~Baker \cite{B}.

Many special cases were considered in the 1980s and early 1990s, and finally in
1996 \spr's conjectures were proved by G.\,A.~Margulis and the author \cite{KM1}
via a method involving dynamics on the space of lattices\footnote{See [Bere] for an
alternative proof and \cite{BKM, BBKM, K} for further developments.}.  In contrast, very
few developments took place in the complex case since \spr's original work. The only
paper on non-polynomial extremal submanifolds of $\bc^n$ known to the author is
\cite{V}, where it is shown that a complex analytic map $\vf:U\to \bc^n$, $U\subset
\bc$, such that 
$$1,f_1,\dots,f_n\text{ are linearly independent
over }\bc\,,\tag 1.5$$ 
is extremal when $n = 3$. Ealier, it was proved in \cite{KoS} that
for arbitrary $n$ and $\vf$ satisfying (1.5), one has $
\omega\big(\vf(z)\big) \le 2n^2 + n - 3$ for almost all $z$.

\smallskip

In the present paper we fill this gap by proving

\proclaim{Theorem 1.1} Let $U\subset\bc^d$ be an open subset, and let $\vf:U\to 
\bc^n$ be a complex analytic map such that {\rm (1.4)} holds. Then $\vf$
is strongly extremal. \endproclaim

\example{Remark 1.2} A `slicing trick' originally due to A.~Pyartli (\cite{P},
see also \cite{S4}) shows that it is enough to prove the above theorem for $d =
1$.
\endexample

\example{Remark 1.3} Note that assumption (1.4) is 
weaker than (1.5). For
example, it follows from Theorem 1.1 that a straight line $\{(z,iz)\mid z\in\bc\}\subset
\bc^2$ is strongly extremal.
\endexample

\example{Remark 1.4} Note also that in the paper \cite{KM1} \spr's conjecture
has been proved in a stronger infinitesimal form, where one  replaces the analyticity 
of $\vf$ by existence of a certain number of derivatives, and (1.4)  by the so-called
`non-degeneracy' condition, equivalent to (1.4) when $\vf$ is real analytic. 
On the other hand, it is clearly impossible to relax the analyticity assumption
in the complex case. For example, $$\vf: z = x + iy \mapsto (z,z^2 \bar z) = (x +
iy)(1,x^2 + y^2)$$  is a map which is polynomial in $x,y$, satisfies (1.5), but is
not extremal.
\endexample

\example{Remark 1.5} In this paper we do not touch the subject of so-called
Khintchine-Groshev-type theorems on complex manifolds, where $\|\vq\|^{-v}$ in the
right hand side of the inequality in (1.2) is replaced by an arbitrary function of
$\|\vq\|$. Results of this type exist for complex polynomials \cite{BD, BernV} and
analytic curves in $\bc^3$ \cite{BereV}. It seems plausible that an approach of this
paper can lead to Khintchine-Groshev-type results for maps $\vf$ as in
Theorem 1.1.
\endexample

\example{Remark 1.6} A separate, and also quite natural, problem is to consider small
values of
$|\vz\cdot\vq + p|$ where both  $p$  and the components of $\vq$ are Gaussian integers,
that is,  integer points of the field
$\bc$. Ore, more generally, one can take two global fields $K\supset L$ with $[K:L] <
\infty$, let
$\bk$ be a completion of $K$ with respect to some valuation (Archimedean or not), and
define 
\de s (multiplicative or not) of 
$\vz\in\bk^n$ with respect to integer points $L_\bz$ of $L$ by looking at small values
of $|\vz\cdot\vq + p|$ where   $p\in L_\bz$  and  $\vq\in (L_\bz)^n$. Note that 
\spr's  book \cite{S3}  solves the analogues of Mahler's Problem for
$L = K$ being
either $\bq$ or the field of formal power series over a finite field. A dynamical
approach
to problems of this type, including the $S$-arithmetic versions where one is allowed to
consider several completions at the same time, is now being developed in
\cite{KT}.
\endexample

Our proof of Theorem 1.1 is based on a variation of methods of \cite{KM1} and
\cite{BKM}. In the next section we describe the reduction of the theorem to a chain of
statements involving discrete subgroups of Euclidean spaces, and in \S 3 take care of the
final  link of that chain.

\heading{2. Lattices and measure estimates}
\endheading 

In order to prove Theorem 1.1, one needs to take any $v > \frac{n-1}2$ and show
that the set of $\vz\in\vf(U)$ for which there exist infinitely many $(p,\vq)\in
\bz^{n+1}$ with
$$
|\vz\cdot\vq + p|   \le
\big(\Pi_+(\vq)\big)^{-v/n}
\tag 2.1
$$
 is null with respect to the pushforward of the Lebesgue measure on $U$.
Writing
$\vz =
\vx + i\vy$ and perhaps slightly changing $v$, one can replace (2.1) with
$$
\max\big(|\vx\cdot\vq + p|,|\vy\cdot\vq| \big)  \le
\big(\Pi_+(\vq)\big)^{-v/n}\,.
\tag 2.2
$$
Our first step is to rephrase (2.2). Choose $\beta > 0$, define
$$
r = \Pi_{\sssize +}(\vq)^{-\beta}\,,\tag 2.3a
$$
 and then define $\vt = (t_1,\dots,t_n)\in
\br^{n}_{\sssize +}$ by 
$$
|q_i|_{\sssize +} = r e^{t_i}\,,\quad i = 1,\dots,n\,.\tag 2.3b
$$
Let us denote the sum of the components of $\vt$ by $t$
(the latter notation will be used throughout the paper, so that whenever 
$t$ and $\vt$ appear in the same formula, $t$ will stand for $\sum_{i = 1}^nt_i$).
Then (2.2a) and (2.2b) imply that
$$
\Pi_{\sssize +}(\vq) = r^ne^{t} = \big(\Pi_+(\vq)\big)^{-n\beta}e^{t} \,,
$$
hence $\Pi_{\sssize +}(\vq) = e^{\frac1{1+n\beta}t}$ and
$$
r = e^{-\frac{\beta}{1+n\beta}t}\,.\tag 2.4
$$
This allows us to write the right hand side of (2.2) as
$$
\big(\Pi_+(\vq)\big)^{-v/n} = e^{-\frac v{1+n\beta}t} =
e^{-\frac{\beta}{1+n\beta}t}e^{-\frac {v-n\beta}{n(1+n\beta)}t} = re^{-at}\,,
$$
where 
$$
a = \frac {v-n\beta}{n(1+n\beta)} \quad \Leftrightarrow 
\quad \beta = \frac {v-an}{n(1+an)}\,.
$$
Now recall that we still have a freedom to choose either $\beta$ or $a$. At this point
we 
choose
$a$ in order to let
$\beta$ 
tend to $0$ as $v$ tends to its critical value $\frac{n-1}2$. That
is, we let $a = \frac{n-1}{2n}$,  which yields
$$
\beta = \frac{2v - n +
1}{n(n+1)}\,,\tag 2.5a
$$
and, in view of (2.4),
$$ r = e^{-\gamma t}\,,\quad \text{where}\quad\gamma =
\frac{2v - n + 1}{2n(v + 1)}\,.\tag 2.5b
$$
We summarize the above computation\footnote{A similar argument can be found in
\cite{KM1, \S 2}, \cite{KM2, \S9}, \cite{K, \S 5}.} as

\proclaim{Lemma 2.1} Let $v > \frac{n-1}2$, $\vx,\vy\in\br^{n}$ and $(p,\vq)\in\bz^{n+1}$
be such that {\rm (2.2)} holds. Define $\beta$ by {\rm (2.5a)}, $r$ by {\rm (2.3a)} and
$\vt$ by {\rm (2.3b)}. Then 
$$
e^{\frac{n-1}{2n}t}\max\big(|\vx\cdot\vq + p|,|\vy\cdot\vq| \big)  \le r\tag 2.6a
$$
 and  
$$
e^{-t_i}|q_i| \le r, \qquad i = 1,\dots,n\,;\tag 2.6b
$$
moreover, $r$ and $t$ are related via {\rm (2.5b)}.
\endproclaim

Now discrete subgroups of $\br^{n+2}$ enter naturally to provide a coincise form for
inequalities (2.6ab). For $\vz = \vx + i\vy\in\bc^n$ define
$$
u_\vz = \left(\matrix
1 & 0 & \vx^{\sssize T}  \\
0 & 1 & \vy^{\sssize T}  \\ 0 & 0 & I_n
\endmatrix \right)\in\SL_{n+2}(\br)\,,\tag 2.7
$$
and for $\vt\in\br^n_{\sssize +}$ let
$$
g_\vt = \text{\rm
diag}(e^{\frac{n-1}{2n}t},e^{\frac{n-1}{2n}t},e^{-t_1},\dots,e^{-t_n})\in\GL_{n+2}(\br)
\tag
2.8
$$
(note that $\det( g_\vt) = e^{-t/n}$). Then (2.6ab) can be rewritten as 
$
\left\|g_\vt u_\vz\left(\matrix
p  \\
0  \\ \vq
\endmatrix \right)\right\|\le r\,,
$
where $\|\cdot \|$ stands for the $l^\infty$ norm on $\br^{n+2}$. Recall
that  for any discrete subgroup $\Lambda$ of $\br^m$, $m\in \bn$, one defines
$\delta(\Lambda)$ to be the norm of a nonzero element of $\Lambda$ with the
smallest norm, that is, 
$$
\delta(\Lambda) \df \inf_{\vv\in\Lambda\nz}\|\vv\|\,.
$$
So if one denotes
$$
\Lambda \df \left\{\left.\left(\matrix
p  \\
0  \\ \vq
\endmatrix \right)\right|p\in\bz,\,\vq\in\bz^n\right\}\,,\tag 2.9
$$
the following is straightforward:

\proclaim{Corollary 2.2}  Let $v > \frac{n-1}2$ and $\vz \in\bc^n$ 
be such that {\rm (2.2)} holds for infinitely many $(p,\vq)\in
\bz^{n+1}$.  Then  there exists an unbounded set of $\vt\in
\br^{n}_{\sssize +}$ 
such that 
$$\delta(g_{\vt}u_\vz\Lambda) \le e^{-\gamma t}\,,\tag 2.10
$$
where $\gamma$ is as in {\rm (2.5b)}.
\endproclaim

\proclaim{Corollary 2.3}  Let  $\vz \in\bc^n$ be such that
for some $v > \frac{n-1}2$ one has {\rm (2.2)}  for infinitely many $(p,\vq)\in
\bz^{n+1}$.  Then  there exists $\gamma > 0$ such that {\rm (2.10)} holds
for infinitely many $\vt\in\bz_{\sssize +}^n$.
\endproclaim

\demo{Proof} Straightforward from (2.5b) and the fact that the ratio of
$\delta(g_{\vt}\cdot)$ and $\delta(g_{\vs}\cdot)$ is uniformly 
bounded from both sides when $\|\vt-\vs\| < 1$.
\qed\enddemo

In fact, the converse to Corollary 2.3 is also true, and can be proved by an argument
from \cite{K, \S 5}.

\smallskip

 In the next corollary and thereafter, $|\cdot|$ stands for Lebesgue
measure.

\proclaim{Corollary 2.4}  Let $\vf$ be a map from an open subset $U$ of $\bc$  to 
$\bc^n$. Suppose that for almost every $z_0\in U$ there exists a neighborhood
$B\subset U$ of $z_0$ such that for any $\gamma > 0$ one has
$$
\sum_{\vt\in\bz^n_{\sssize +}}\big|\big\{z\in B\bigm|
\delta\big(g_{\vt}u_{\vf(z)}\Lambda\big)
\le e^{-\gamma t}\big\}\big| < \infty\,.\tag 2.11
$$
 Then
$\vf$ is strongly extremal. \endproclaim

\demo{Proof} In view of  the Borel-Cantelli Lemma, it follows from (2.11) 
that for almost all $\vz$ of the form $\vf(z)$, $z\in B$,  (2.10) is satisfied
for at most finitely many $ \vt\in\bz_{\sssize +}^n$. Corollary 2.3 then
implies that for any $v > \frac{n-1}2$, almost all $\vz$ as above
satisfy (2.1) for at most finitely many $(p,\vq)\in
\bz^{n+1}$, that is, $\big|\big\{z\in B\bigm|
\omega^\times\big(\vf(z)\big) \ge v\big\}\big| = 0$.
\qed\enddemo

In view of the above corollary and Remark 1.2, to prove Theorem 1.1 it suffices
to show that for any complex analytic $\vf:U\to 
\bc^n$, $U\subset\bc$, satisfying (1.4), one can find a neighborhood
$B\subset U$ of almost every $z_0\in U$  such that  for any $\gamma > 0$  and 
any $\vt\in
\bz^{n}_{\sssize +}$ it is possible to estimate the measure of sets
$$\{z\in B\bigm|
\delta\big(g_{\vt}u_{\vf(z)}\Lambda\big)
\le e^{-\gamma t}\big\}\tag 2.12
$$
so that   (2.11) holds. Observe that at this point it makes no difference if
one replaces the $l^\infty$ norm used to define $\delta(\cdot)$ by any other norm,
and we are going to switch to the Euclidean one from now on.

\smallskip

In order to 
state the main estimate that will be  used to bound  the measure of sets
(2.12), we need to introduce some additional terminology.  If
$\Gamma$ is a discrete subgroup of
$\br^m$, we define the {\sl rank\/} of $\Gamma$ to be the dimension of  $\br\Gamma$, 
and denote
by $\|\Gamma\|$ the {\sl covolume\/} of $\Gamma$, that is, the volume of the quotient
space
$\br\Gamma/\Gamma$. For
$\Lambda$ as above, we  denote by $\Cal S(\Lambda)$ the set of all nonzero subgroups of
$\Lambda$. The key ingredient in what follows is, for fixed $\vt\in\bz^n_{\sssize +}$ and
$\Gamma\in\Cal S(\Lambda)$,  keeping track of the covolumes of subgroups
$g_{\vt}u_{\vf(z)}\Gamma$ as functions of
$z$. In particular, it will be useful to see that all those functions share
a certain property, referred to in \cite{KM1} and subsequent papers as being \cag.

Namely, if  $C$ and $\alpha$ are  positive
numbers and $V$ is a subset of $\br^d$, one says that   a function 
$f:V\mapsto \br$ is {\sl \cag\ on\/}   $V$  if for any ball
$B\subset V$ and any $\vre > 0$
one has 
$$
\big|\{\vx\in B\bigm| |f(\vx)| < \vre\cdot{\sup_{\vx\in B}|f(\vx)|}\}\big| \le
C\vre^\alpha |B|\,. 
$$

Later we will need the following facts:

 \proclaim{Lemma 2.5} 
{\rm (a)} Suppose that   $f_1,\dots,f_k$
are $(C,\alpha)$-good on
$V$; then the function $(f_1^2 + \dots +f_k^2)^{1/2}$ is $(k^{\alpha/2}C,\alpha)$-good
on
$V$.

{\rm (b)} Let $\vf$ be a
real analytic map from a connected open subset $U$ of
$\br^d$ to $\br^n$. Then for any  $\vx_0\in
U$ there exists a neighborhood
$V \subset U$ of $\vx_0$ and positive $C,\alpha$  such that any linear combination of $
1,f_1,\dots,f_n$ is
$(C,\alpha)$-good on
$V$. 
 \endproclaim  

\demo{Proof} The first statement is elementary and an easy consequence of parts (b), (c)
of \cite{BKM, Lemma 3.1}. For the second assertion, which is based on the work done in
\cite{KM1}, see \cite{K, Corollary 3.3}. \qed
\enddemo

Now we can state the crucial estimate that the whole proof hinges upon. It is 
a special case of \cite{BKM, Theorem 6.2}. We remark that it is proved by a
variation of a combinatorial construction used by Margulis in the 1970s
\cite{Ma} and then by S.\,G.~Dani in the 1980s \cite{D} to establish and quantitatively
describe non-divergence of unipotent flows on \hs s.

 \proclaim{Theorem 2.6} Fix $\,d,k,m\in\bn$. Let  $\Lambda$ be a discrete subgroup of
$\br^m$ of rank $\,k$,  and let a ball $B = B(\vx_0,r_0)\subset \br^d$ and a continuous
map
$H:\tilde B
\to \GL_m(\br)$ be given, where $\tilde B$ stands for $B(\vx_0,3^kr_0)$. Take
$C,\alpha > 
0$, $0 < \rho  
\le 1$, and assume that for any $\Gamma\in\Cal S(\Lambda)$,

{\rm(i)} the function $\vx\mapsto \|H(\vx)\Gamma\|$ is \cag\ on
$\tilde B$,  and  

{\rm(ii)} 
$\sup_{\vx\in B}\|H(\vx)\Gamma\| \ge \rho$.

\noindent Then 
 for any  positive $ \vre \le \rho$ one has 
$$
\left|\big\{\vx\in B\bigm| \delta\big(H(\vx)\Lambda\big) < \vre  \big\}\right| \le
c_{d,k} C \left(\frac\vre
\rho \right)^\alpha  |B|\,,
$$
where $c_{d,k}$ is a constant depending only on $d$ and $k$ (and explicitly computed in
\cite{KM1} and \cite{BKM}). 
\endproclaim

\proclaim{Corollary 2.7}   Let
$U\subset\bc$ be an open subset, and let
$\vf:U\to 
\bc^n$ be a continuous  map. Keep the notation {\rm (2.7), (2.8), (2.9)}. Assume that
for almost every $z_0\in U$ one can find  balls
$B= B(z_0,r_0)\subset \tilde B = B(z_0,3^{n+1}r_0)\subset U$ and constants $C,\alpha,
\rho >  0$ such
that for any $\vt\in\br^n_{\sssize +}$ and any $\Gamma\in\Cal S(\Lambda)$ the following
holds:
$$
\text{$\forall\,\vt\in\br^n_{\sssize +}\ \forall\,\Gamma\in\Cal S(\Lambda)\quad$ 
the function $z\mapsto \|g_{\vt}u_{\vf(z)}\Gamma\|$ is \cag\ on
$\tilde B$}\,,\tag 2.13
$$
and
$$
\text{$\forall\,\vt\in\br^n_{\sssize +}\ \forall\,\Gamma\in\Cal S(\Lambda)\quad
\sup_{z\in B}\|g_{\vt}u_{\vf(z)}\Gamma\| \ge \rho$}\,.\tag 2.14
$$

\noindent  Then
$\vf$ is strongly extremal. \endproclaim

\demo{Proof} Under the assumptions (2.13) and (2.14) above, the previous theorem, with
$d = 2$, $k = n+1$, $m = n+2$ and $H(z) = g_{\vt}u_{\vf(z)}$ for fixed
$\vt\in\bz^n_{\sssize +}$, forces  the measure of every set (2.12) to be  not greater
than $c_{d,k} C\rho^{-\alpha}  |B| e^{-\alpha\gamma t}$ whenever $\vt$ is far enough
from zero.
This implies  (2.11), and hence, in view of Corollary 2.4, the strong extremality of
$\vf$.
\qed\enddemo

In the next section 
we show how one can write down explicit expressions for the covolume of  subgroups of the
form
$g_{\vt}u_{\vf(z)}\Gamma$ for any $\Gamma\in\Cal S(\Lambda)$, and, assuming (1.4) and the
analyticity of $\vf$,  verify conditions (2.13) and (2.14). 


\heading{3. Exterior products and covolume estimates}
\endheading

In this section we   keep the notation introduced in {\rm (2.7)--(2.9)},
and prove the following:

\proclaim{Lemma 3.1}   Let
$U\subset\bc$ be an open subset, and let
$\vf:U\to 
\bc^n$ be a complex analytic map. Then:

{\rm (a)} for any $z\in U$ there exists a ball $\tilde B\subset U$ centered at $z$ and
$C,\alpha > 0$ such that {\rm (2.13)} holds;

{\rm (b)} if, in addition, {\rm (1.4)} is satisfied, then for every ball $B\subset U$
there exists $\rho > 0$ such that {\rm (2.14)} holds.
 \endproclaim

In view of Corollary 2.7 and Remark 1.2, this lemma immediately implies Theorem 1.1.

\demo{Proof} 
 Fix $\vt\in\br^n_{\sssize +}$ and
$\Gamma\in\Cal S(\Lambda)$ of rank $k$, where  $1\le k \le n+1$. Without loss of
generality we will order the components of
$\vt$ so that
$t_1 \le\dots\le t_n$.

It will be convenient to denote the 
standard basis of  
$\br^{n+2}$ by $\{\ve_0,\ve_*,\ve_1,\dots,\ve_n\}$, so that $\Lambda$ as in
(2.9) is equal to $\bz\ve_0 + \bz\ve_1 + \dots +\bz\ve_n$. We will describe $\Gamma$ by
means of its  representing element from  the exterior algebra of
$\br^{n+2}$. Recall that 
 $\vw\in \bigwedge^k(\br^{n+2})$ is said to {\sl represent\/} $\Gamma$ if 
$$
\vw =  \vv_{1}\wedge\dots\wedge \vv_{k}\,, \text{ where }\vv_{1},\dots, \vv_{k} \text{
form a basis of }
\Gamma\,.
$$
Clearly an element representing $\Gamma$ is defined up to a sign, and the
covolume $\|\Gamma\|$ of $\Gamma$ is 
equal to
the  Euclidean norm (with respect to the standard Euclidean structure  extended to
$\bigwedge(\br^{n+2})$) of $\vw$.
 
\smallskip

For brevity we will suppress the variable $z$
whenever it does not cause confusion. Write $\vf = \vg + i\vh$, and with some
abuse of notation identify
$\vg$ and
$\vh$ with vector-functions 
$\pmatrix 0\\0\\
\vg\endpmatrix$ and $\pmatrix 0\\0\\
\vh\endpmatrix: U\to\br^{n+2}$.  Then it is
immediate from (2.7) that for any
$\vv\in\br^{n+2}$ one has
$$u_\vf\vv =
 \vv + (\vv\cdot \vg)\ve_0 + (\vv\cdot \vh)\ve_*\,.\tag 3.1
$$

Our goal now is to choose an orthonormal set
with respect to which it is 
convenient to compute coordinates of 
$u_\vf\vw$, where $\vw$ represents
$\Gamma$. We closely follow an approach developed in \cite{BKM}, where statements
similar to (2.13) and (2.14) were established  to
prove 
the convergence case of the \kgt\ for non-degenerate
manifolds.

\smallskip

First choose an
orthonormal subset
$\vv_1,\dots,\vv_{k-1}$ of
$\br\Gamma$ such that  each $\vv_i$, $i = 1,\dots,k-1$, is orthogonal to $\ve_0$. Then,
if $\br\Gamma$ does not contain $\ve_0$, choose 
$\vv_0\in \br\ve_0\oplus \br\Gamma$ such that 
$\{\ve_0,\vv_0, \vv_1,\dots,\vv_{k-1}\}$ is  an orthonormal  basis of $\br\ve_0\oplus
\br\Gamma$, and  represent $\Gamma$ by 
$$
\vw = (a \ve_{0}  + b\vv_{0})\wedge\vv_{1}\dots\wedge \vv_{k-1} =  a
\ve_{0}\wedge\vv_{1}\dots\wedge \vv_{k-1} + b
\vv_{0}\wedge\vv_{1}\wedge\dots\wedge \vv_{k-1}\,,\tag 3.2 
$$
where $a^2 + b^2 \ge 1$. If $\br\Gamma$ does  contain $\ve_0$, then
$\{\ve_0,\vv_1,\dots,\vv_{k-1}\}$ is already a basis of $\br\ve_0 + \br\Gamma$, so
(3.2) is valid with $b = 0$ and   $\vv_0$ taken to be any unit vector
orthogonal to $
\br\ve_*\oplus\br\Gamma$. 

Combining  (3.1) and (3.2), one gets 
$$
u_\vf\vw   =  \Big(a \ve_{0}  + b\big(\vv_{0} + (\vv_0\cdot \vg)\ve_0 + (\vv_0\cdot
\vh)\ve_*\big)\Big) \wedge \bigwedge_{i = 1}^{k-1}\big(\vv_{i} + (\vv_i\cdot \vg)\ve_0 +
(\vv_i\cdot
\vh)\ve_*\big)\,.\tag 3.3
$$
From there one can see that every coordinate of $
u_\vf\vw$ with respect to the basis
$$
\left\{\bigwedge_{i = 0} ^{k-1}  
\vv_i,\ \ve_0\wedge\bigwedge_{s \ne
i}  
\vv_s,\ \ve_*\wedge\bigwedge_{s \ne
i}  
\vv_s, \ \ve_0\wedge\ve_*\wedge\bigwedge_{s \ne
i,j}  
\vv_s\right\}\tag 3.4
$$
of $\bigwedge^k(\br\ve_0\oplus
\br\ve_* + \br\Gamma + \br\vv_0)$ is a linear combination of functions
$$
1,\ \vv_i\cdot \vg, \ \vv_i\cdot \vh, \ \left|\matrix \vv_i \cdot
\vg & \vv_j\cdot  \vg\\ \vv_i 
\cdot \vh & \vv_j \cdot \vh\endmatrix \right|\,,
$$
hence a linear combination of $1$ and the components of the map $\tilde
\vf:U\to\br^{n(n+3)/2}$ given  by
$$
\tilde \vf\df 
\left(\vg,\ \vh\,;\ \left|\matrix g_i & g_j\\ h_i 
 & h_j \endmatrix \right|\,,\ 1\le i < j \le n\right)\,.
$$
The same can be said about every coordinate of $
g_\vt u_\vf\vw$ with respect to any orthonormal  basis of $\bigwedge^k(\br^{n+2})$
containing the one given by  (3.4). Since all the components of $\tilde \vf$ are real
analytic functions, it follows from Lemma 2.5(b) that for any $z\in U$ there exists a
ball 
$\tilde B\subset U$ centered at $z$ and
$C',\alpha > 0$ such that every linear combination of $1$ and the components of $\tilde
\vf$ is $(C',\alpha)$-good on $\tilde B$. Part (a) of this lemma then immediately
follows from Lemma
2.5(a).

\smallskip

Part (b) requires some more work. Let us state the following
auxiliary result:

 \proclaim{Lemma 3.2} Let $B\subset \bc$ be a nonempty ball.

{\rm (a)} Suppose that   $\vf = (f_1,\dots,f_n)$ is an $n$-tuple of real-valued 
functions
on
$B$ satisfying {\rm (1.4)}. Then there  exists $\rho_1 = \rho_1(B,\vf) > 0$ such that for
any 
$\vc =( c_0,c_1,
\dots , c_n)$ with $\|\vc\| = 1$ one has
$$\sup_{z\in
B}|c_0 + c_1f_1(z) + \dots + c_nf_n(z)| \ge \rho_1\,.
$$

{\rm (b)} Suppose that $\Cal F$ is a compact (in $C^0$ topology) family of pairs
of complex-valued functions $(\varphi_1,\varphi_2)$  analytic in $B$ 
such that one is not a real multiple of another. 
Then there  exists $\rho_2 = \rho_2(B,\Cal F) > 0$ such that
$$|\im(\bar{\varphi}_1\varphi_2)|\ge
\rho_2 \quad \forall\,(\varphi_1,\varphi_2)\in \Cal F\,.
$$
 \endproclaim  

\demo{Proof} The first statement follows from the  standard compactness argument applied
to the set of functions $\big\{c_0 + c_1f_1 + \dots + c_nf_n\bigm| \|\vc\| =
1\big\}$. The same kind of argument shows that if  the second assertion does not hold,
one must have $\im(\bar{\varphi}_1\varphi_2) = 0$ for some $(\varphi_1,\varphi_2)\in
\Cal F$. Hence the ratio of $\varphi_1(z)$ and $\varphi_2(z)$ is real for all $z\in B$,
which, due to the complex analyticity of $\varphi_1$ and $\varphi_2$, can happen only if 
 $\varphi_1/\varphi_2$ is a constant.
\qed
\enddemo

Now  let us consider two cases. If $k = \text{dim}(\br\Gamma) = 1$, equality (3.3)
gives
$$
u_\vf\vw   =  \big(a   + b(\vv_0\cdot \vg)\big)\ve_0  + b\vv_{0} +  b(\vv_0\cdot
\vh)\ve_* \,,
$$
hence 
$$
\|g_{\vt}u_\vf\vw \| \ge |g_{\vt}u_\vf\vw\cdot\ve_0| = e^{\frac{n-1}{2n}t}|a   +
b(\vv_0\cdot \vg)|\,.\tag 3.5
$$
The right hand side of (3.5) is a linear combination of $1,g_1,\dots,g_n$ (which, due to
(1.4), are linear independent over $\br$) with big enough coefficients, therefore its
supremum on any ball $B\subset \bc$ 
is uniformly (in $a,b,\vv_0$ and $\vt$) bounded from below by $\rho_1(B,\vg)$ due to
Lemma 3.2(a).

\smallskip

The argument in the case $k\ge 1$ is different. Namely, we define the family
$$
\Cal F \df \big\{\big(\vu_1\cdot \vf, a + b\vu_2 \cdot \vf\big)\mid a^2 + b^2 = 1,\ 
\vu_1 \perp \vu_2 \in \br^n,\  \|\vu_1\| = \|\vu_2\| = 1  \big\}\,,\tag 3.6
$$
which clearly  satisfies all the assumptions of Lemma 3.2(b), and put $\rho =
\rho_2(B,\Cal F)$. 

 To prove that (2.14) holds with this value of $\rho$, we need to fine-tune
the choice of  the orthonormal set
$\{\vv_0,\dots,\vv_{k-1}\}$. Namely, we will pay special  attention to the vector
$\ve_n$, which is the eigenvector of $g_\vt$ with one of the smallest eigenvalues
(recall that $t_1 \le\dots\le t_n$ by our assumption). We do it by first choosing an
orthonormal set
$\vv_1,\dots,\vv_{k-2}\in \br\Gamma$ such that  each $\vv_i$, $i = 1,\dots,k-2$, is
orthogonal to both $\ve_0$ and $\ve_n$. Then choose $\vv_{k-1}$ orthogonal to $\vv_i$,
$i = 1,\dots,k-2$, and to $\ve_0$ (but in general not to $\ve_n$). After that choose
$\vv_0$  either to
complete
$\{\ve_0,\vv_1,\dots,\vv_{k-1}\}$ to an orthonormal basis of $\br\ve_0\oplus
\br\Gamma$ or (in case $\br\Gamma$ contains $\ve_0$) to be any unit vector
orthogonal to $
\br\ve_*\oplus\br\Gamma$.

With this in mind,  denote by $W$ the subspace $\ve_0  \wedge
\ve_*\wedge\bigwedge^{k-2}(\br^{n+2})$ of $\bigwedge^{k}(\br^{n+2})$
(that is, the set of elements corresponding to $k$-dimensional subspaces
of $\br^{n+2}$ containing both $\ve_0$ and
$\ve_*$), and write (3.3) in
the form
$$
u_{\vf(z)}\vw   =  \vw'(z) + \ve_0  \wedge \ve_*\wedge\vw''(z)\,,
$$
where  $\vw'(z)$ is orthogonal to $W$  and 
$\vw''(z)\in\bigwedge^{k-2}\big((\br\ve_0
\oplus
\br\ve_*)^\perp\big)$. Since both
$\ve_{0}$ and $\ve_*$ are eigenvectors of $g_\vt$,
both $W$ and its orthogonal complement  are
$g_\vt$-invariant.  Therefore it will suffice  to show that
$$
\sup_{z\in B}\|g_\vt\big(\ve_0  \wedge \ve_*\wedge\vw''(z)\big)\|\ge\rho\,,
$$
or, equivalently, $$
\sup_{z\in B}\|g_\vt\vw''(z)\|\ge e^{-\frac{n-1}{n}t}\rho\,.
$$

Next consider the product $\ve_n\wedge\vw''(z) $. We claim that it is enough to show
that  
$$
\|\ve_n\wedge\vw''(z) \| \ge \rho\quad\text{for some }z\in B\,.\tag 3.7
$$
 Indeed, since 
 $\ve_n$ is an eigenvector of $g_\vt$ with eigenvalue $e^{-t_n}$, for any $z\in B$ the
norm of $g_\vt\big(\ve_n\wedge\vw''(z)\big)$ is not greater than
$e^{-t_n}\|g_\vt\vw''(z)\|$. Therefore, since the smallest eigenvalue of $g_\vt$ on
$\bigwedge^{k-1}(\br^{n+2})$ is equal to $e^{-(t_{n-k+1} +
\dots +t_n)}$ (here we set $t_0 = -\frac{n-1}{2n}t$, so that the above statement holds
for $k = n+1$ as well), the norm of
$g_\vt\vw''(z)$ is not less than 
$$
\split
&e^{t_n}\|g_\vt\big(\ve_n\wedge\vw''(z)\big)\| 
\ge e^{t_n}e^{-(t_{n-k+1} +
\dots +t_n)}\|\ve_n\wedge\vw''(z)\| \\
\un{\text{for some $z\in B$, by (3.7)}}
\ge \ &e^{t_n}e^{-(t_{n-k+1} +
\dots +t_n)}{\rho}\ 
\ge  \ e^{t_n-t}{\rho} \ge e^{-\frac{n-1}{n}t}\rho
\,,
\endsplit
$$
as required.
Thus it remains to prove (3.7).  Using (3.3), it is possible to write down coefficients
in the decomposition of 
 $\vw''(z)$ as a linear combination of elements of the form $\bigwedge_{s \ne
i,j}  
\vv_s$. We are going to do it for  the term containing
$\vv_1\wedge\dots\wedge \vv_{k-2}$. Namely, one can write 
$$
\aligned
\ve_n\wedge\vw'' &=  \pm \left|\matrix \vv_{k-1}\cdot 
\vg & a + b\vv_0 \cdot \vg \\  \vv_{k-1} \cdot \vh & b\vv_0 
\cdot \vh \endmatrix \right|\,\ve_n\wedge \vv_1\wedge\dots\wedge \vv_{k-2} \\
+ \text{ other terms }&\text{where one or two of $\vv_i$, $i = 1,\dots,k-2$, are
missing.}
 \endaligned \tag 3.8
$$
From the orthogonality of the two summands in (3.8)
it follows that $\|\ve_n\wedge\vw''\|$ is not less than the norm of the
first summand, which, in view of $\ve_n$ being orthogonal to $\vv_i$, $i =
1,\dots,k-2$, is equal to the absolute value of the coefficient in front of 
$\ve_n\wedge \vv_1\wedge\dots\wedge \vv_{k-2}$. The latter can be written as
$$\sqrt{a^2 + b^2}\im(\bar{\varphi}_1\varphi_2)$$ where 
$${\varphi}_1 = \vv_{k-1}\cdot 
\vf\quad\text{and}\quad{\varphi}_2 = (a^2 + b^2)^{-1/2}(a + b \vv_{0}\cdot 
\vf)\,.$$ 
It is immediate that $(\varphi_1,\varphi_2)$ belongs to $\Cal F$ as defined in (3.6).
Since $a^2 + b^2 \ge 1$, this completes the proof of (3.7) with 
$\rho =
\rho_2(B,\Cal F)$. \qed
\enddemo

\heading{Acknowledgements}
\endheading 

The author is  grateful to G.\,A.~Margulis for  many stimulating  discussions, and
to the hospitality of ETH/Zurich where a substantial part of this work has been
accomplished.

\enddocument